# FUZZY LOGIC CONTROL OF A HYBRID ENERGY STORAGE MODULE FOR NAVAL PULSED POWER APPLICATIONS

Isaac Cohen[1], David Wetz[1], Stepfanie Veiga[2], Qing Dong[2], and John Heinzel[2]

[1]Electrical Engineering Department, University of Texas at Arlington, Arlington, USA
[2]Naval Surface Warfare Center, Philadelphia Division, Philadelphia, USA

## ABSTRACT

*There is need for an energy storage device capable of transferring high power in transient situations aboard naval vessels. Currently, batteries are used to accomplish this task, but previous research has shown that when utilized at high power rates, these devices deteriorate over time causing a loss in lifespan. It has been shown that a hybrid energy storage configuration is capable of meeting such a demand while reducing the strain placed on individual components. While designing a custom converter capable of controlling the power to and from a battery would be ideal for this application, it can be costly to develop when compared to purchasing commercially available products. Commercially available products offer limited controllability in exchange for their proven performance and lower cost point - often times only allowing a system level control input without any way to interface with low level controls that are frequently used in controller design. This paper proposes the use of fuzzy logic control in order to provide a system level control to the converters responsible for limiting power to and from the battery. A system will be described mathematically, modeled in MATLAB/Simulink, and a fuzzy logic controller will be compared with a typical controller.*

## 1. INTRODUCTION

The development of shipboard electrical systems continues to grow as the Navy plans to utilize more electrical components than previously seen. With this growth comes a new electrical load profile not commonly seen aboard ships. The addition of new directed energy based weapons as well as electric propulsion give the ship a load profile that is more transient in nature than would typically be seen. Traditionally, power generation sources would be fossil fuel driven generators, but such a configuration would result in very poor performance and power quality if used to drive the transient load profiles envisioned for the future [1].

Previous research has shown that energy storage devices (ESDs) are capable at compensating the stochastic and intermittent nature of these power demands by absorbing the excessive energy when generation exceeds predicted levels and providing energy back to the power system network when generation levels fail to meet the demand. These ESDs are an essential component to future power systems when integrating variable energy resources and stochastic pulsed loads. Previous work has shown that Hybrid Energy Storage Modules (HESMs) can contribute to not only improve the performance of an ESD, but also overcome the limitations of the individual components of the architecture, such as the power delivery limitations of a battery or the energy storage limitations of an ultracapacitor [2, 3, 4]. While this topology has been verified, there are still some questions on how to best control them. [5] designed a low cost digital energy management system and an optimal control algorithm was developed in [6] to coordinate slow ESDs and fast ESDs. While these control developments have been very useful, they do not address the need for a controller that is capable of interfacing with Commercial-Off-The-Shelf (COTS) products. Fuzzy logic control has been used in many applications such as [7], who developed an intuitive fuzzy logic based learning algorithm which was implemented to reduce intensive computation of a complex dynamics such as a humanoid. Some, such as [8] have utilized fuzzy logic control to drive a HESM, but in their case, they used the controller in order to eliminate the need to constantly calculate resource intensive Riccati equations to assist in choosing gains for an adaptive Linear-Quadratic Regulator controller. Others, such as [9, 10] developed an energy-based split and power sharing control strategy for hybrid energy storage systems, but these strategies are focused on different target variables such as loss reduction, leveling the components state of charge, or optimizing system operating points in a vehicular system.

When constructing a HESM for naval pulsed power load applications, several design parameters must be considered to develop the controller. First, the HESM must be responsible for supplying all load demands that it might encounter. Energy sizing techniques such as shown in [11] should be applied in order to assist in meeting this demand. Second, the HESM should add the benefit of becoming an additional source during parallel operation with an existing shipboard power system, offering a reduction in external power system sizing. Finally, the HESM should operate as a power buffer during this parallel operation, offering improved power quality to the load from sluggish generation sources, such as fossil-fuel generators, such as seen in [12].

There are many combinations of high power density and high energy density ESDs that can be used to create a HESM, but a common topology that will be evaluated in this paper is the combination of ultracapacitors and batteries. Although the Navy has traditionally used lead-acid batteries as the go-to ESD, lithium-ion batteries (LIBs) offer a higher energy and power density with respect to both volume and weight. In this scenario, a HESM becomes even more necessary as previous research has shown that while LIBs are capable of being cycled at high rates, they degrade much more quickly [13] – necessitating the limiting the power to and from the battery to be as minimal as possible. This limitation can be implemented through the introduction of intelligent controls.

While this paper aims to address the need for a system level control for COTS devices, in order to simulate the controller in MATLAB/Simulink, it is necessary to use a simplified model of a HESM and its power converters. Other researchers have spent time modeling this system mathematically such as [14], which presented a detailed small-signal mathematical model that represents the dynamics of the converter interfaced energy storage system around a steady-state operating point. Their model considered the variations in the battery current, supercapacitor current, and DC load bus voltage as the state variables, the variations in the power converter's duty cycle as the input, and the variations in the battery voltage, supercapacitor voltage, and load current as external disturbances. This paper will utilize MATLAB/Simulink's SimPowerSystems toolbox instead of a mathematical model and will be explained further in Section II.

The paper is organized as follows: Section II will provide a system description for the HESM model being used in MATLAB/Simulink's SimPowerSystems and the fuzzy logic controller designed to accomplish the aforementioned tasks. Section III will present the simulation's results and Section IV will provide concluding remarks.

## 2. SYSTEM DESCRIPTION

### 2.1. Modeling the HESM

A generic schematic of a HESM is shown in Figure 1. In this schematic, a simple buck-boost converter is utilized in order to give the controller a method of bi-directional voltage and current control. The load and the generator are tied together as one variable current source/sink as the system which the HESM augments can be seen as a generalized external power disturbance. When mathematically modeling this system, in similar fashion to [14], C1 and C2, the output capacitors for each direction of power flow, hold the state variables of the battery bus voltage and the DC load bus voltage and L, the power inductor, holds the state variable of the power converter current, and the combined current sourcing or sinking from the load and generation is the external disturbance, denoted as $\zeta$. For instance, when $\zeta>0$, the generator is producing more current than the load is drawing, but when $\zeta<0$, the load is drawing more current than the generation is producing.

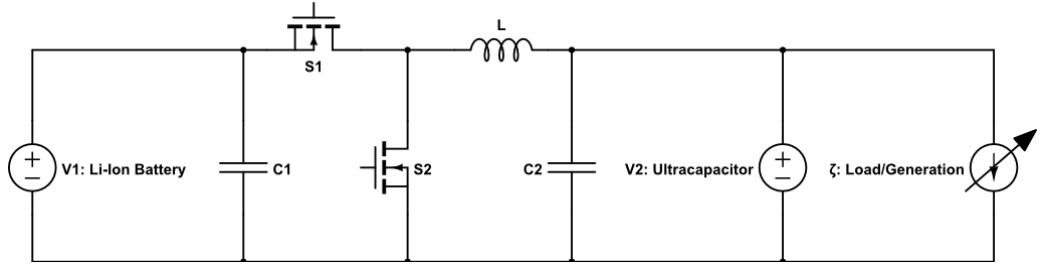

Figure 1: Schematic of a generic battery/ultracapacitor HESM

There are multiple equations to describe this circuit during operation, based on the state of the switches. To demonstrate the variation of this system over time, the mathematical equations that represent these states are shown below. For simplicity, the system will be evaluated in each direction of power flow while treating the load ESD as an omitted independent variable, shown in the circuit as $V_1$ and $V_2$.

1. State 1: When power flows from the battery towards the ultracapacitor and when $S_1$ is on and $S_2$ is off,

$$L\frac{di_L}{dt} = V_{C1} - V_{C2} \qquad (1)$$

$$C_1\frac{dV_{C1}}{dt} = i_L - i_{V1} \qquad (2)$$

$$C_2\frac{dV_{C2}}{dt} = i_L - i_{V2} - \zeta \qquad (3)$$

2. State 2: When power flows from the battery towards the ultracapacitor and when $S_1$ is off and $S_2$ is off,

$$L\frac{di_L}{dt} = -V_{C2} \qquad (4)$$

$$C_1\frac{dV_{C1}}{dt} = -i_{V1} \qquad (5)$$

$$C_2\frac{dV_{C2}}{dt} = i_L - i_{V2} - \zeta \qquad (6)$$

3. *State 3: When power flows from the ultracapacitor towards the battery and when $S_1$ is off and $S_2$ is on,*

$$L\frac{di_L}{dt} = V_{C2} \tag{7}$$

$$C_1\frac{dV_{C1}}{dt} = -i_{V1} \tag{8}$$

$$C_2\frac{dV_{C2}}{dt} = i_{V2} - i_L - \zeta \tag{9}$$

4. *State 4: When power flows from the ultracapacitor towards the battery and when $S_1$ is off and $S_2$ is off,*

$$L\frac{di_L}{dt} = V_{C2} - V_{C1} \tag{10}$$

$$C_1\frac{dV_{C1}}{dt} = i_L - i_{V1} \tag{11}$$

$$C_2\frac{dV_{C2}}{dt} = i_{V2} - i_L - \zeta \tag{12}$$

These mathematical equations are presented here to enforce the idea that this system is time varying, with four different plant descriptions, depending on the state of 2 switches. Although there are mathematical methods in which these equations could be combined to produce a time-average model in which the controller could be evaluated, it was decided that a more accurate model could be created using MATLAB/Simulink's SimPowerSystems toolbox for controller evaluation in order to utilize the toolbox's lithium ion battery model and additional component values such as internal impedances. The block diagram for this model can be seen in Figure 2.

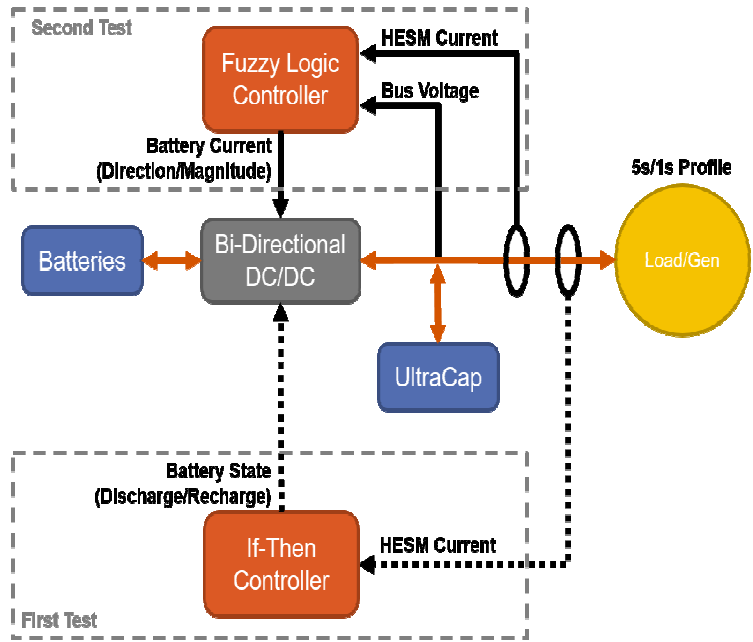

Figure 2: Block diagram of HESM experiment in SimPowerSystems toolbox

Initial investigation into acquiring parts for hardware validations lead to the values chosen for this simulation. This investigation occurred with the intention of validating this controller with real hardware for future investigations. The values chosen for the simulation can be seen in Table 1.

Table 1: Simulation component values

| Component | Value |
|---|---|
| *Lithium Ion Battery* | 36 $V_{nom}$, 15 Ahr, 50% SOC |
| *Ultracapacitor* | 24 $V_{nom}$, 29 F, 44 mΩ ESR |
| *Load/Generation* | 5 seconds on, 1 second off, variable currents |
| *Switches* | $R_{on}$ = 100 mΩ, $f_{sw}$ = 40 kHz |
| *Inductor* | 3.4 mH, 1.5 mΩ ESR |
| *Discharge PI Controller* | Voltage: $k_P$ = 0.2, $k_I$ = 10<br>Current: $k_P$ = 1, $k_I$ = 200 |
| *Recharge PI Controller* | Voltage: $k_P$ = 0.028, $k_I$ = 1.5<br>Current: $k_P$ = 5, $k_I$ = 1 |

## 2.2. Designing the Fuzzy Logic Controller

When designing a HESM, one of the largest obstacles to overcome is the successful implementation of system level control. Fuzzy Logic Control (FLC) employs an if-then rule-base with mathematical fuzzification and defuzzification in order to achieve an expert response with a digital controller's speed and efficiency. Fuzzy systems typically achieve utility in assessing more conventional and less complex systems [15], but on occasion, FLC can be useful in a situation where highly complex systems only need approximated and rapid solutions for practical applications. FLC can be particularly useful in nonlinear systems such as this HESM which shifts between 4 different operation states. One key difference between crisp and fuzzy sets is their membership functions. The uniqueness of a crisp set is sacrificed for the flexibility of a fuzzy set. Fuzzy membership functions can be adjusted to maximize the utility for a particular design application. The membership function embodies the mathematical representation of membership in a set using notation $\Omega_i$, where the functional mapping is given by $\mu_{\Omega_i}(x) \in [0,1]$. The symbol $\Omega_i(x)$ is the degree of membership of element x in fuzzy set $\Omega_i$ and $\mu_{\Omega_i}(x)$ is a value on the unit interval which measures the degree to which $x$ belongs to fuzzy set $\Omega_i$.

Two fuzzy input sets were defined as the DC bus voltage and the HESM current – these inputs can be seen in Figure 3 and 4. The first input, the DC bus voltage, was a logical choice as maintaining this voltage is critical to all three tasks of the HESM, which is to supply power under all scenarios and to act as a power buffer to transients, as transients will cause deviations in the DC bus voltage. The second input, HESM current, was chosen because of its proportional relativity to the differential change in the bus voltage. The HESM current can be defined as,

$$i_{HESM} = i_{C2} + i_L \qquad (13)$$

which leads to,

$$\frac{i_{HESM} - i_L}{C_2} = \frac{dV_{C2}}{dt} \qquad (14)$$

where $V_{C2}$ is the DC bus voltage. It is because of this proportional relationship that the HESM current is able to give the controller sensory insight to the direction of the demand for power, giving an increased ability to maintain the voltage of the DC bus, similar to what a traditional PID controller would be able to offer.

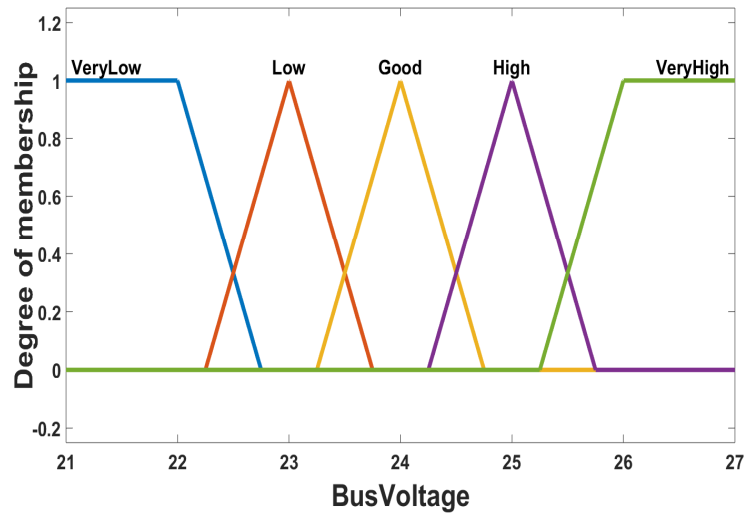

Figure 3: First input fuzzy membership functions

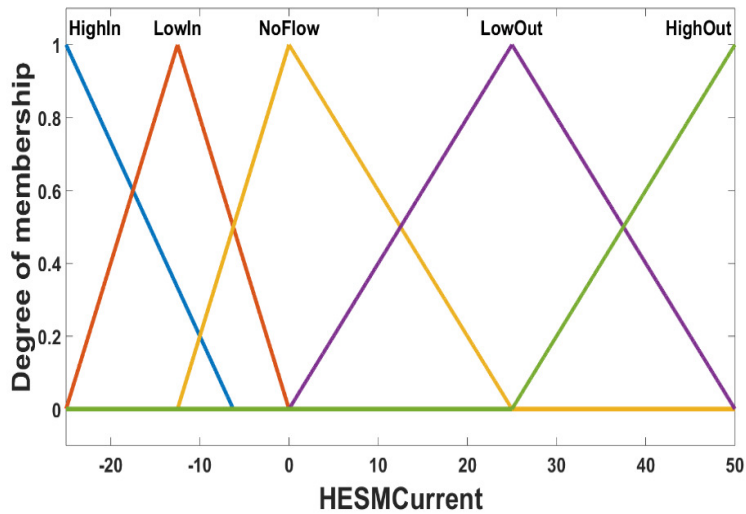

Figure 4: Second input fuzzy membership functions

The output fuzzy set was chosen to be the battery current limit as this gives a superb amount of control over the HESM. The membership function describing this relationship can be seen in Figure 5. With the battery being actively limited in the direction and magnitude, all remaining HESM power must come from the ultracapacitor. This essentially gives full, albeit indirect, control over the power flow to and from all ESDs associated with the HESM.

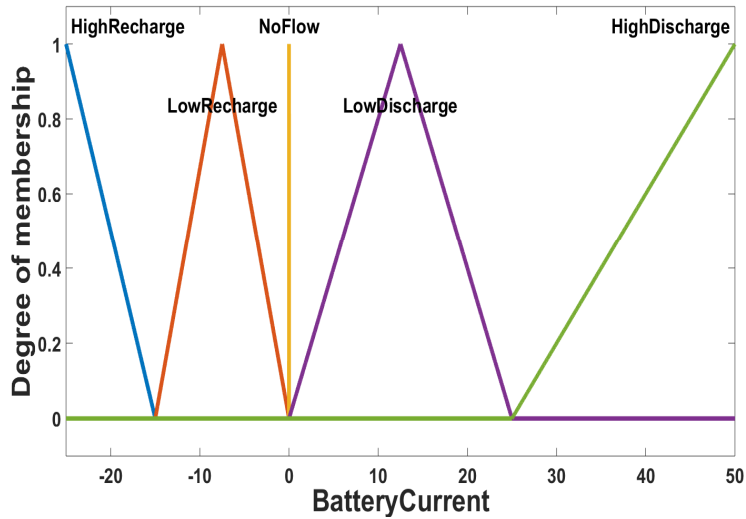

Figure 5: Output fuzzy membership function

The controller's inputs and outputs are mapped together using a rule-base. That is to say that for a given set of input conditions, there should be a relative output condition. Since it is possible for values to lie in-between conditions, there may be times when multiple output conditions are met. In this case, the FLC determines the value by computing the centroid of mass in the membership functions [16, 17]. The set of rules being used in this controller is seen in Table 2. The process of relating these input values to a rule-base is called fuzzification.

Table 2: Fuzzy Logic Rule-Base

|  |  | Bus Voltage | | | | |
| --- | --- | --- | --- | --- | --- | --- |
|  |  | Very Low | Low | Good | High | Very High |
| HESM Current | High Out | No Flow | Low Recharge | High Recharge | High Recharge | High Recharge |
|  | Low Out | Low Discharge | No Flow | Low Recharge | High Recharge | High Recharge |
|  | No Flow | High Discharge | Low Discharge | No Flow | Low Recharge | High Recharge |
|  | Low In | High Discharge | High Discharge | Low Discharge | No Flow | Low Recharge |
|  | High In | High Discharge | High Discharge | High Discharge | Low Discharge | No Flow |

Once the rule-base is in place, the controller's input versus output can be observed as seen in Figure 6. The plateau area represents high discharge while the valley area represents high recharge. The process of relating the output membership function to a value is called defuzzification.

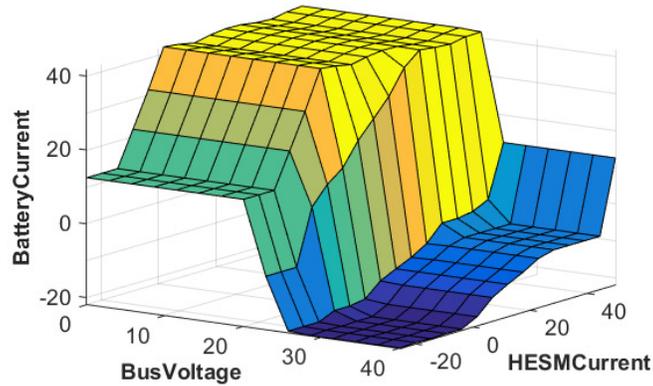

Figure 6: Fuzzy Logic Controller input vs. output relationship

### 2.3. Experiment Description

In order to test the controller an experiment was designed to mimic a typical naval pulsed power load. In this case, a pulse train load profile of 5 seconds at high load power draw and 1 second of low load power draw was simulated. To compare the FLC to a traditional type of state-machine controller, the experiment introduces a shift in the power demand halfway through the test.

## 3. SIMULATION RESULTS

For comparison, Figure 7 shows the results of the experiment when using a simple if-then controller and Figure 8 shows the results of the experiment when using the FLC. Positive currents indicate that a device is sourcing energy and negative currents indicate that a device is sinking energy.

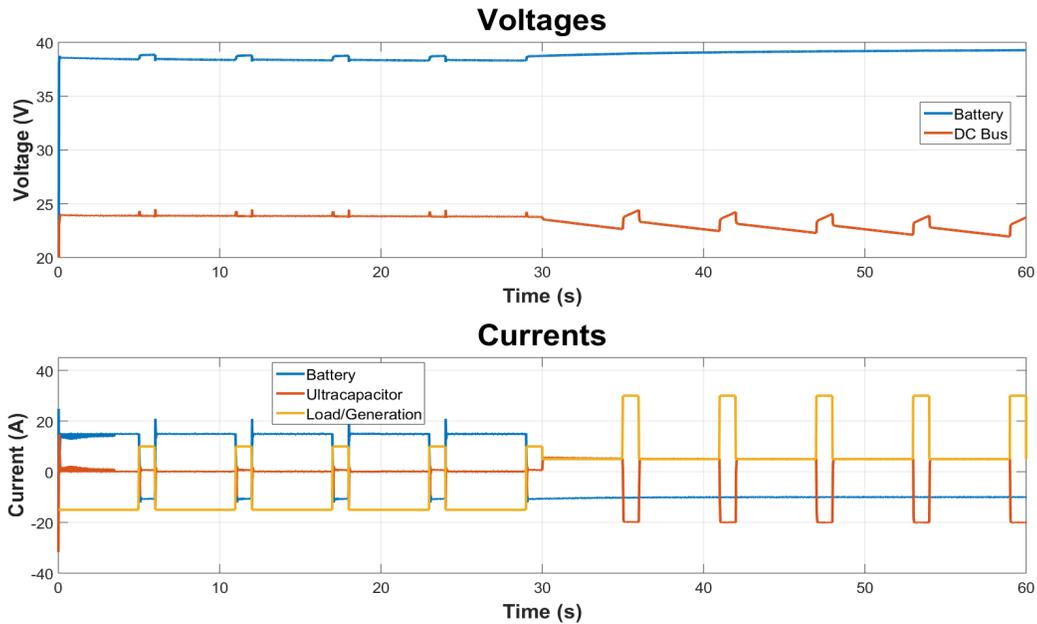

Figure 7: Simulation results with if-then control

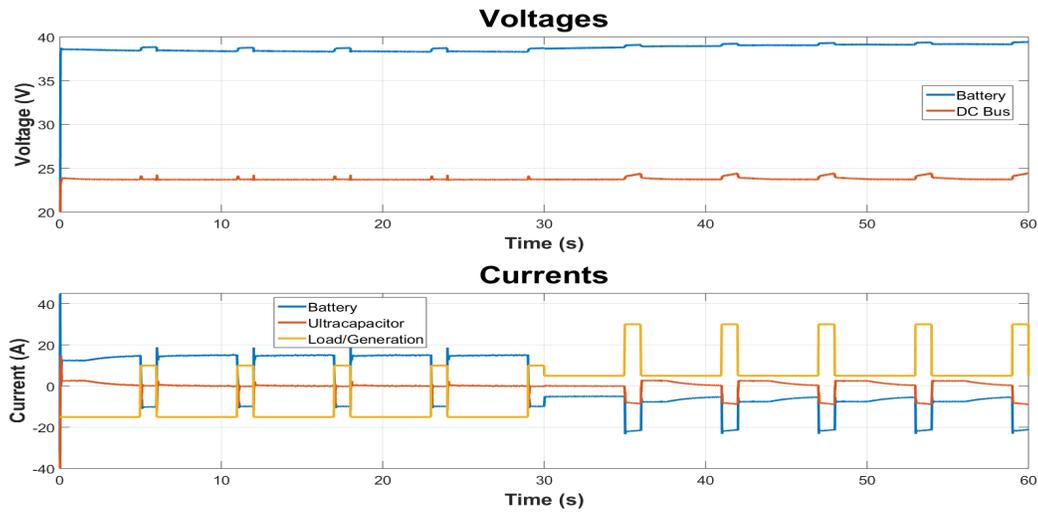

Figure 8: Simulation results using Fuzzy Logic Control

By examination of Figure 7, it is clear to see that the controller operates satisfactory for the first 30 seconds, where it was designed to operate, but as soon as the system starts to exhibit behaviors outside of a pre-determined need, the controller is no longer able to effectively maintain the DC bus voltage. In addition to the presence of large voltage swings during the change in the pulsed load power draw, there is an overall decay in voltage, which in a real system would lead to a failure to maintain an effective power buffer and therefore lead to either a cascading power failure throughout other components of the integrated power system from over demand, or a larger sizing requirement to accommodate loads outside of pre-determined profiles.

Comparing the results in Figure 7 to Figure 8, it is clear to see that as the load/generation shifts to a different region, the FLC is able to easily accommodate the change. There are certainly some changes in the voltage swing through the pulsed power load profile, but this can be fixed through more meticulous tuning of the FLC. A numerical representation of the results can be seen in Table 3.

Table 3: Numerical experimental results - comparison of if-then control with FLC

| | | |
|---|---|---|
| **If-Then** | *Voltage Swing during Pulse* | 1.77 V |
| | *Voltage Sag over final 30 seconds* | 1.84 V |
| **FLC** | *Voltage Swing during Pulse* | 0.66 V (62% improvement) |
| | *Voltage Sag over final 30 seconds* | 0.00 V (100% improvement) |

## 4. CONCLUSIONS

ESDs are becoming more and more crucial as integrated power systems evolve. Their application in naval settings are becoming more desirable and the challenges associated with individual ESDs can be overcome by utilizing a HESM topology. This paper proposed one method of system level control in this topology and has demonstrated through simulation that this is a viable technique. One drawback of the FLC is that it is only as accurate as the knowledge of the expert creating it, but through meticulous tuning, even this inadequacy can be overcome.

## ACKNOWLEDGEMENTS

This material is based upon work supported by the US Office of Naval Research (ONR). The authors would like to express thanks to ONR for their continued support. Any opinions, findings, and conclusions or recommendations expressed in this publication are those of the authors and do not necessarily reflect the views of the US Office of Naval Research.

# Authors

**Isaac J. Cohen** was born in Miami, FL in 1988. He received the B.S. degree in electrical engineering in 2013 from the University of Texas at Arlington, where he is currently working toward the Ph.D. degree in electrical engineering.He served as Chair of the Student Branch of IEEE at UTA in 2012, where he received several awards for his outstanding service. He currently serves as Chair of the Young Professional affinity group within the IEEE Fort Worth Section.His research interests include applying control theories to power electronics in microgrid, energy storage, and pulsed power settings. 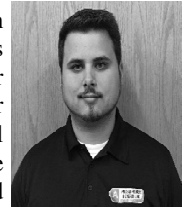

**David A. Wetz, Jr.** was born in El Paso, TX, USA, in 1982. He received the B.Sc. degrees in electrical engineering and computer science, the M.Sc. degree in electrical engineering, and the Ph.D. degree in electrical engineering from Texas Tech University, Lubbock, TX, USA, in 2003, 2004, and 2006, respectively.He was a Post-Doctoral Fellow with the Institute for Advanced Technology, University of Texas at Austin, Austin, TX, USA, from 2006 to 2007, where he was also a Research Associate from 2007 to 2010. He joined the Faculty of Electrical Engineering at the University of Texas at Arlington, Arlington, TX, USA, as an Assistant Professor in 2010. His current research interests include pulsed power, power electronics, energy storage, and power system analysis. 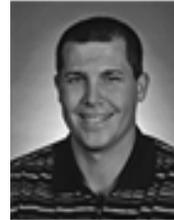

**Stepfanie Veiga** Received the M.S. degree in Electrical Engineering from Villanova University, Villanova, PA, USA, in 2013. She spent 4 years in private industry at firms such as Specialty Minerals Inc., Bethlehem, PA, USA, and Bard, Rao + Athanas Consulting Engineers, LLC, Philadelphia, PA, USA. She is currently with the Machinery Technology Research & Development Group, Philadelphia Division, Naval Surface Warfare Center, U.S. Navy, Philadelphia, PA, USA, where she has been for 7 years. She also has prior experience in Shipboard Instrumentation and System Calibration for the Machinery Information, Sensors, and Control Systems division at Naval Surface Warfare Center, Philadelphia. Her current research interests include power electronics, converter topologies, pulsed power system, high voltage power systems testing and evaluation, and real-time power hardware in the loop architecture development and testing. 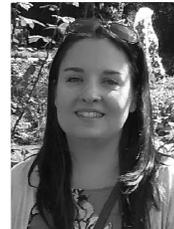

**Qing Dong** received the Ph.D. degree in electrical engineering from Temple University, Philadelphia, PA, USA, in 2011. His thesis topic was titled Multi-Agent Based Federated Control of Large-Scale Systems. He spent 16 years in private industry at firms such as Lucent Technologies, Murray Hill, NJ, USA, and Parker Hannifin Corporation, Cleveland, OH, USA. He is currently with the Machinery Technology Research & Development Group, Philadelphia Division, Naval Surface Warfare Center, U.S. Navy, Philadelphia, PA, USA, where he has been for 12 years. He has over a dozen peer-reviewed publications in journals and conference proceedings. His current research interests include optimal control, large-scale system dynamics and control, agent-based control methods, and optimization theory. 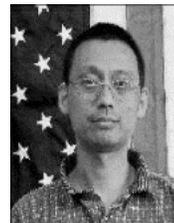

**John M. Heinzel** received the bachelor's degree from the University of Delaware, Newark, DE, USA, and the Ph.D. degree in chemical engineering from Auburn University, Auburn, AL, USA. He is currently a Senior Chemical Engineer with the Carderock Division, Naval Surface Warfare Center, U.S. Navy, Philadelphia, PA, USA. He also serves as the lead on multiple Navy/DoD programs in these areas, and emphasizes a Multiphysics approach toward solving problems and transitioning new hardware into practical applications. His area of cognizance is chemical systems related to energy conversion research and development, in particular, advanced fuel cell and energy storage systems, and other power generation, directed energy, and thermal management areas. 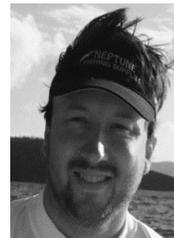